\documentclass[reqno,12pt,a4paper]{amsart}

\usepackage{mathrsfs}
\usepackage{amssymb}
\usepackage{amsfonts}
\usepackage{latexsym}
\usepackage{amsthm}

\usepackage{graphicx}
\def\lb{\label}

\newcommand{\er}[1]{\textrm{(\ref{#1})}}

\begin{document}


\renewcommand{\theequation}{\arabic{section}.\arabic{equation}}
\theoremstyle{plain}
\newtheorem{theorem}{\bf Theorem}[section]
\newtheorem{lemma}[theorem]{\bf Lemma}
\newtheorem{corollary}[theorem]{\bf Corollary}
\newtheorem{proposition}[theorem]{\bf Proposition}
\newtheorem{definition}[theorem]{\bf Definition}
\newtheorem{remark}[theorem]{\it Remark}

\def\a{\alpha}  \def\cA{{\mathcal A}}     \def\bA{{\bf A}}  \def\mA{{\mathscr A}}
\def\b{\beta}   \def\cB{{\mathcal B}}     \def\bB{{\bf B}}  \def\mB{{\mathscr B}}
\def\g{\gamma}  \def\cC{{\mathcal C}}     \def\bC{{\bf C}}  \def\mC{{\mathscr C}}
\def\G{\Gamma}  \def\cD{{\mathcal D}}     \def\bD{{\bf D}}  \def\mD{{\mathscr D}}
\def\d{\delta}  \def\cE{{\mathcal E}}     \def\bE{{\bf E}}  \def\mE{{\mathscr E}}
\def\D{\Delta}  \def\cF{{\mathcal F}}     \def\bF{{\bf F}}  \def\mF{{\mathscr F}}
\def\c{\chi}    \def\cG{{\mathcal G}}     \def\bG{{\bf G}}  \def\mG{{\mathscr G}}
\def\z{\zeta}   \def\cH{{\mathcal H}}     \def\bH{{\bf H}}  \def\mH{{\mathscr H}}
\def\e{\eta}    \def\cI{{\mathcal I}}     \def\bI{{\bf I}}  \def\mI{{\mathscr I}}
\def\p{\psi}    \def\cJ{{\mathcal J}}     \def\bJ{{\bf J}}  \def\mJ{{\mathscr J}}
\def\vT{\Theta} \def\cK{{\mathcal K}}     \def\bK{{\bf K}}  \def\mK{{\mathscr K}}
\def\k{\kappa}  \def\cL{{\mathcal L}}     \def\bL{{\bf L}}  \def\mL{{\mathscr L}}
\def\l{\lambda} \def\cM{{\mathcal M}}     \def\bM{{\bf M}}  \def\mM{{\mathscr M}}
\def\L{\Lambda} \def\cN{{\mathcal N}}     \def\bN{{\bf N}}  \def\mN{{\mathscr N}}
\def\m{\mu}     \def\cO{{\mathcal O}}     \def\bO{{\bf O}}  \def\mO{{\mathscr O}}
\def\n{\nu}     \def\cP{{\mathcal P}}     \def\bP{{\bf P}}  \def\mP{{\mathscr P}}
\def\r{\rho}    \def\cQ{{\mathcal Q}}     \def\bQ{{\bf Q}}  \def\mQ{{\mathscr Q}}
\def\s{\sigma}  \def\cR{{\mathcal R}}     \def\bR{{\bf R}}  \def\mR{{\mathscr R}}
\def\S{\Sigma}  \def\cS{{\mathcal S}}     \def\bS{{\bf S}}  \def\mS{{\mathscr S}}
\def\t{\tau}    \def\cT{{\mathcal T}}     \def\bT{{\bf T}}  \def\mT{{\mathscr T}}
\def\f{\phi}    \def\cU{{\mathcal U}}     \def\bU{{\bf U}}  \def\mU{{\mathscr U}}
\def\F{\Phi}    \def\cV{{\mathcal V}}     \def\bV{{\bf V}}  \def\mV{{\mathscr V}}
\def\P{\Psi}    \def\cW{{\mathcal W}}     \def\bW{{\bf W}}  \def\mW{{\mathscr W}}
\def\o{\omega}  \def\cX{{\mathcal X}}     \def\bX{{\bf X}}  \def\mX{{\mathscr X}}
\def\x{\xi}     \def\cY{{\mathcal Y}}     \def\bY{{\bf Y}}  \def\mY{{\mathscr Y}}
\def\X{\Xi}     \def\cZ{{\mathcal Z}}     \def\bZ{{\bf Z}}  \def\mZ{{\mathscr Z}}
\def\O{\Omega}
\def\ve{\varepsilon}   \def\vt{\vartheta}    \def\vp{\varphi}
\def\vk{\varkappa}     \def\vs{\varsigma}

\newcommand{\mc}{\mathscr {c}}

\newcommand{\gA}{\mathfrak{A}}          \newcommand{\ga}{\mathfrak{a}}
\newcommand{\gB}{\mathfrak{B}}          \newcommand{\gb}{\mathfrak{b}}
\newcommand{\gC}{\mathfrak{C}}          \newcommand{\gc}{\mathfrak{c}}
\newcommand{\gD}{\mathfrak{D}}          \newcommand{\gd}{\mathfrak{d}}
\newcommand{\gE}{\mathfrak{E}}
\newcommand{\gF}{\mathfrak{F}}           \newcommand{\gf}{\mathfrak{f}}
\newcommand{\gG}{\mathfrak{G}}           
\newcommand{\gH}{\mathfrak{H}}           \newcommand{\gh}{\mathfrak{h}}
\newcommand{\gI}{\mathfrak{I}}           \newcommand{\gi}{\mathfrak{i}}
\newcommand{\gJ}{\mathfrak{J}}           \newcommand{\gj}{\mathfrak{j}}
\newcommand{\gK}{\mathfrak{K}}            \newcommand{\gk}{\mathfrak{k}}
\newcommand{\gL}{\mathfrak{L}}            \newcommand{\gl}{\mathfrak{l}}
\newcommand{\gM}{\mathfrak{M}}            \newcommand{\gm}{\mathfrak{m}}
\newcommand{\gN}{\mathfrak{N}}            \newcommand{\gn}{\mathfrak{n}}
\newcommand{\gO}{\mathfrak{O}}
\newcommand{\gP}{\mathfrak{P}}             \newcommand{\gp}{\mathfrak{p}}
\newcommand{\gQ}{\mathfrak{Q}}             \newcommand{\gq}{\mathfrak{q}}
\newcommand{\gR}{\mathfrak{R}}             \newcommand{\gr}{\mathfrak{r}}
\newcommand{\gS}{\mathfrak{S}}              \newcommand{\gs}{\mathfrak{s}}
\newcommand{\gT}{\mathfrak{T}}             \newcommand{\gt}{\mathfrak{t}}
\newcommand{\gU}{\mathfrak{U}}             \newcommand{\gu}{\mathfrak{u}}
\newcommand{\gV}{\mathfrak{V}}             \newcommand{\gv}{\mathfrak{v}}
\newcommand{\gW}{\mathfrak{W}}             \newcommand{\gw}{\mathfrak{w}}
\newcommand{\gX}{\mathfrak{X}}               \newcommand{\gx}{\mathfrak{x}}
\newcommand{\gY}{\mathfrak{Y}}              \newcommand{\gy}{\mathfrak{y}}
\newcommand{\gZ}{\mathfrak{Z}}             \newcommand{\gz}{\mathfrak{z}}

\def\A{{\mathbb A}} \def\B{{\mathbb B}} \def\C{{\mathbb C}}
\def\dD{{\mathbb D}} \def\E{{\mathbb E}} \def\dF{{\mathbb F}} \def\dG{{\mathbb G}}
\def\H{{\mathbb H}}\def\I{{\mathbb I}} \def\J{{\mathbb J}} \def\K{{\mathbb K}}
\def\dL{{\mathbb L}}\def\M{{\mathbb M}} \def\N{{\mathbb N}} \def\dO{{\mathbb O}}
\def\dP{{\mathbb P}} \def\R{{\mathbb R}}\def\S{{\mathbb S}} \def\T{{\mathbb T}}
\def\U{{\mathbb U}} \def\V{{\mathbb V}}\def\W{{\mathbb W}} \def\X{{\mathbb X}}
\def\Y{{\mathbb Y}} \def\Z{{\mathbb Z}}
\def\dk{{\Bbbk}}


\def\la{\leftarrow}              \def\ra{\rightarrow}            \def\Ra{\Rightarrow}
\def\ua{\uparrow}                \def\da{\downarrow}
\def\lra{\leftrightarrow}        \def\Lra{\Leftrightarrow}


\def\lt{\biggl}                  \def\rt{\biggr}
\def\ol{\overline}               \def\wt{\widetilde}
\def\no{\noindent}


\let\ge\geqslant                 \let\le\leqslant
\def\lan{\langle}                \def\ran{\rangle}
\def\/{\over}                    \def\iy{\infty}
\def\sm{\setminus}               \def\es{\emptyset}
\def\ss{\subset}                 \def\ts{\times}
\def\pa{\partial}                \def\os{\oplus}
\def\om{\ominus}                 \def\ev{\equiv}
\def\iint{\int\!\!\!\int}        \def\iintt{\mathop{\int\!\!\int\!\!\dots\!\!\int}\limits}
\def\el2{\ell^{\,2}}             \def\1{1\!\!1}
\def\sh{\sharp}
\def\wh{\widehat}
\def\bs{\backslash}
\def\intl{\int\limits}

\def\na{\mathop{\mathrm{\nabla}}\nolimits}
\def\sh{\mathop{\mathrm{sh}}\nolimits}
\def\ch{\mathop{\mathrm{ch}}\nolimits}
\def\where{\mathop{\mathrm{where}}\nolimits}
\def\all{\mathop{\mathrm{all}}\nolimits}
\def\as{\mathop{\mathrm{as}}\nolimits}
\def\Area{\mathop{\mathrm{Area}}\nolimits}
\def\arg{\mathop{\mathrm{arg}}\nolimits}
\def\const{\mathop{\mathrm{const}}\nolimits}
\def\det{\mathop{\mathrm{det}}\nolimits}
\def\diag{\mathop{\mathrm{diag}}\nolimits}
\def\diam{\mathop{\mathrm{diam}}\nolimits}
\def\dim{\mathop{\mathrm{dim}}\nolimits}
\def\dist{\mathop{\mathrm{dist}}\nolimits}
\def\Im{\mathop{\mathrm{Im}}\nolimits}
\def\Iso{\mathop{\mathrm{Iso}}\nolimits}
\def\Ker{\mathop{\mathrm{Ker}}\nolimits}
\def\Lip{\mathop{\mathrm{Lip}}\nolimits}
\def\rank{\mathop{\mathrm{rank}}\limits}
\def\Ran{\mathop{\mathrm{Ran}}\nolimits}
\def\Re{\mathop{\mathrm{Re}}\nolimits}
\def\Res{\mathop{\mathrm{Res}}\nolimits}
\def\res{\mathop{\mathrm{res}}\limits}
\def\sign{\mathop{\mathrm{sign}}\nolimits}
\def\span{\mathop{\mathrm{span}}\nolimits}
\def\supp{\mathop{\mathrm{supp}}\nolimits}
\def\Tr{\mathop{\mathrm{Tr}}\nolimits}
\def\BBox{\hspace{1mm}\vrule height6pt width5.5pt depth0pt \hspace{6pt}}


\newcommand\nh[2]{\widehat{#1}\vphantom{#1}^{(#2)}}
\def\dia{\diamond}

\def\Oplus{\bigoplus\nolimits}



\def\qqq{\qquad}
\def\qq{\quad}
\let\ge\geqslant
\let\le\leqslant
\let\geq\geqslant
\let\leq\leqslant
\newcommand{\ca}{\begin{cases}}
\newcommand{\ac}{\end{cases}}
\newcommand{\ma}{\begin{pmatrix}}
\newcommand{\am}{\end{pmatrix}}
\renewcommand{\[}{\begin{equation}}
\renewcommand{\]}{\end{equation}}
\def\eq{\begin{equation}}
\def\qe{\end{equation}}
\def\[{\begin{equation}}
\def\bu{\bullet}

\title[{1D Schr\"odinger operators with complex potentials }]
{1D Schr\"odinger operators with complex potentials}

\date{\today}

\author[Evgeny Korotyaev]{Evgeny Korotyaev}
\address{Saint-Petersburg State University,
Universitetskaya nab. 7/9, St. Petersburg, 199034, Russia, \
E-mail address:
korotyaev@gmail.com, \ e.korotyaev@spbu.ru }


\begin{abstract}
\no  We consider a Schr\"odinger operator with complex-valued
potentials  on the line. The operator  has essential spectrum on the
half-line plus eigenvalues (counted with algebraic multiplicity) in
the complex plane without the positive half-line. We determine
series of trace formulas. Here we have the new term: a singular
measure, which is absent for real potentials. Moreover, we estimate
of sum of Im part of eigenvalues plus singular measure in terms of
the norm of potentials. The proof is based on classical results
about the Hardy spaces.

\end{abstract}

\maketitle


\section {Introduction and main results}
\setcounter{equation}{0}

\subsection{Introduction}
We consider a Schr\"odinger operator $H=-{d^2\/dx^2}+q(x)$ on the
space $L^2(\R)$, where the potential $q$  is complex and satisfies:
\[
\lb{dV} \int_\R(1+|x|)|q(x)|dx<\iy.
\]
 It is known that the spectrum of
the operator $H$ has two components: the essential spectrum which covers the
half-line $[0,\iy)$
 plus $N\le \iy$ eigenvalues (counted with multiplicity) in the cut
 spectral domain $\C\sm [0,\iy)$. We denote them by
 $
 E_j\in \C\sm [0,\iy), j=1,...,N,
 $
according to their multiplicity. Note, that the multiplicity of each
eigenvalue equals 1, but we call the multiplicity of the eigenvalue
its algebraic multiplicity. Define the half-planes $\C_\pm=\{\pm \Im
z>0\}$. Instead of the energy $E\in \C$ we define the momentum
$k=\sqrt E\in \ol\C_+$. We call $k_j=\sqrt E_j\in \C_+$ also the
eigenvalues of the operator $H$. Of course, $E$ is really the
energy, but since $k$ is the natural parameter, we will abuse
terminology. We  define the set $\dk_q=\{k_1,...,k_{N}\in \C_+\}$
and label $ k_1,..,k_{N}\in \C_+$ by
\[
\lb{Lkn} \Im k_1\ge \Im k_2\ge \Im k_3\ge... \ge \Im k_n\ge ...
\]

We shortly describe results about  trace formulas:

$\bu $ In 1960 Buslaev and Faddeev  \cite{BF60} determined the
classical results about trace formulas for Schr\"odinger operators
with real decaying potentials on half-line. The case of the real
line was discussed by Faddeev and Zakharov in the nice paper
\cite{FZ71}.

$\bu $ There are a lot of results about one dimensional case, see
\cite{KS09} and references therein.

$\bu $ The  multidimensional case was studied in \cite{B66}. Trace
formulas for Stark operators and magnetic Schr\"odinger operators
were discussed in \cite{KP03}, \cite{KP04}.

$\bu $ The trace formulas for Schr\"odinger operators with real
periodic potentials were obtained in \cite{KK95,K97}. They were used
to obtain two-sided estimates of potential in terms of gap lengths
(or the action variables for KdV) in \cite{K00} via the conformal
mapping theory for the quasimomentum.

$\bu $   Trace formulas for Schr\"odinger operators with complex
potentials on the lattice $\Z^d$ and on $\R^3$ are considered
recently in \cite{K17}, \cite{KL18}, \cite{MN15}   and  \cite{K17x}
respectively. In \cite{K17}, \cite{KL18} for  the discrete case  the
main tool is the classical results about the Hardy spaces in the
disc. Trace formulas for Schr\"odinger operators with complex
potentials on  $\R_+$  with the Dirichlet boundary condition are
discussed in \cite{K18}. In the case $\R^3$ and $\R_+$ the Hardy
spaces in the upper half-plane $\C_+$ are used.

\subsection{The Hardy spaces}
  Introduce the Jost solutions $f_\pm(x,k )$  of the equation
\[
\lb{1.3} -{f_\pm}''+qf_\pm=k ^2f_\pm,\ \ \ x\in \R, \ \ \ k \in
\ol\C_+\sm \{0\} ,
\]
with the conditions
\[
\begin{aligned}
\lb{pf} f_\pm(x,k)=e^{\pm ixk }+o(1)\qqq as \qq x\to \pm\iy, \qqq
k\in \R\sm \{0\}.
\end{aligned}
\]
Here and in the following $'$ denotes the derivative w.r.t. the
first variable. For each $x\in \R$ the Jost solutions $f^\pm(x,k)$
are analytic in $\C_+ $, continuous up to the real line. Introduce
the Wronskian $w$ and functions $\p, \P$ in $\C_+$ by
\[
w(k)=\{f_-(x,k), f_+(x,k)\}|_{x=0},\qqq \p(k)={w(k)\/2ik},\qq
\P={w\/2i(k+i)},
\]
where $\{y, f\}=yf'-y'f$. The function $\p$ satisfies
(uniformly in $\arg k \in [0,\pi]$):
\[
\lb{asa1} \p(k)=1-{q_0+o(1)\/2ik}\qqq \as \qqq |k |\to \iy,\qq
q_0:=\int_\R q(t)dt.
\]
The function $\p$ has $N\geq 0$ zeros in $\C_+$ given by $k_{j}\in
\dk_q$, counted with multiplicity.

\medskip

Define the Hardy space $\mH_p$. Let a function $F(k), k=u+iv\in
\C_+$ be analytic on $\C_+$. For $0<p\le \iy$ we say $F\in
\mH_p=\mH_p(C_+)$ if $F$ satisfies $\|F\|_{\mH_p}<\iy$, where
$\|F\|_{\mH_p}$ is given by
$$
\|F\|_{\mH_p}=\ca
\sup_{v>0}{1\/2\pi}\rt(\int_\R|F(u+iv))|^pdu\rt)^{1\/p} &
if  \qqq 0< p<\iy\\
 \sup_{k\in \C_+}|F(k)| & if \qqq p=\iy\ac .
$$
Note that the definition of the Hardy space $\mH_p$ involves all
$v=\Im k>0$.

We define
$$
 \|q\|_\a=\int_\R  |x|^\a |q(x)|dx, \ \a\ge 0,\qqq \|q\|=\|q\|_0.
$$

Describe the properties of the functions $w,\p, \P$:

\no $\bu $  The function $w,\p, \P$ have the same zeros in $\ol
\C_+\sm \{0\}$,  and these zeros $\{k_j\}$ in the upper-half plane
$\C_+$ labeled by \er{Lkn} satisfy (see e.g., \cite{Sa10}):
\[
\lb{B1} \sum _{j=1}^N \Im k_j<\iy.
\]
\no  $\bu $ if $w(0)\ne 0$ and $q_0\ne0$, then   the function
$\p-1\notin \mH_p$ for any $p>0$, since $\p$ has asymptotics
\er{asa1}.

\no $\bu $ The properties of $w$ give that $\P-1\in \mH_p$ for any
$p>1$ and all $q$.

\no $\bu $ If $w(0)=0$, then   the function $\p-1\in \mH_p$ for any
$p>1$.

\no $\bu $ Due to  \er{asa1} all zeros of $w$ are uniformly bounded
and satisfy (see e.g., \cite{AAD01})
\[
\lb{1.2}
\begin{aligned}
\dk_q\ss \{k\in \ol\C_+: w(k)=0\}     \ss \{k\in \ol\C_+: |k|\le
r_{c}\},\qq  \where\qq  r_{c}:={\|q\|\/2}.
 \end{aligned}
\]
In order to study $\p(k)$ in the upper-half plane we defined the
Blaschke product $B$  by
\begin{equation}
\label{Bk}
 B(k)=\prod_{j=1}^N\rt(\frac{k-k_j}{ k-\ol k_j}\rt),\qqq k\in \C_+.
\end{equation}
This product converges absolutely for each $k\in \C_+$, since all
zeros of $w$ are uniformly bounded, see \er{1.2}. Moreover, it has
an analytic continuation from $\C_+$ into the domain $\{|k|>r_c\}$,
where $r_c={\|q\|\/2}$ and has the following Taylor series
\[
\begin{aligned}
\lb{B6} & \log  B(k)=-i{B_0\/k}-i{B_1\/2k^2}-i{B_2\/3k^3}-..., \qqq
as \qqq |k|>r_c,
\\
 & B_0=2\sum_{j=1}^N\Im k_j,\qqq B_n=2\sum_{j=1}^N\Im k_j^{n+1},\qq n\ge 1,
\end{aligned}
\]
where each sum $B_n, n\ge 1$ is absolutely convergence and satisfies
\[
|B_n|\le 2\sum_{j=1}^N|\Im k_j^{n+1}|\le {\pi}(n+1)r_c^{n}B_0.
\]
We use asymptotics of $B$ at large $|k|$  to determine the trace
formulas.  Note that the function $B$ has a complicated properties
in the disk $\{|k|<r_c\}$ and very good properties for
$\{|k|>r_c\}$.
 We describe the basic properties of eigenvalues and the Blaschke product
$B$.

\begin{proposition}\lb{T1}
Let a  potential $q$ be complex and  satisfy
$\int_\R(1+|x|)|q(x)|dx<\iy$. Then

i) The Blaschke product $B(k), k\in \C_+$ given by \er{Bk} belongs
to $\mH_\iy$ with $\|B\|_{\mH_\iy}\le 1$.

ii) Let $q_0=\int_\R q(x)dx, \ \|q\|_1=\int_\R |xq(x)|dx $ and let
$A=\|q\|\|q\|_1e^{\|q\|_1}$. Then

$\bu$ If $A<\Re q_0$, then the operator $H$ does  not have
eigenvalues.

$\bu$  If $A<\Re (-q_0)$, then the operator $H$ has exactly one
simple eigenvalue.

\end{proposition}

\no {\bf Example.} Consider the potential $q(x)=ct x^{t^2-1}, x\in
(0,1)$ and $q(x)=0$ for $x>1$, where $c\in \C, t>0$. We have
$$
\|q\|={|c|\/t},\qqq \|q\|_1={|c|t\/1+t^2},\qqq q_0={c\/t},\qqq
A={|c|^2\/1+t^2}\exp {|c|t\/1+t^2}.
$$
If $t$ is small, then the complex potential $q$ is rather big and
due to \er{1.2} all eigenvalues belong to the half-disk with the
radius $r_c={\|q\|\/2}={|c|\/2t}$. Let $c=|c|e^{i\f}$ and let $t>0$
be sufficiently small. We have two cases:

If $\f={\pi\/3}$, then $\Re q_0={|c|\/2t}>A$
 and by Proposition \ref{T1}, the operator $H$ has not
eigenvalues.

If $\f=-{\pi\/3}$, then $\Re (-q_0)={|c|\/2t}>A$  and by Proposition
\ref{T1}, the operator $H$ has one simple eigenvalue.

\subsection{Trace formulas and estimates}
We describe the function $\p$ in terms of a canonical factorization.
In general, the function $\p\notin \mH_p$ for all $p>0$, but we show
that  the function $\p$ has  a canonical factorization for each
potential $q$.

\begin{theorem}
\label{T2} Let a potential $q$ satisfy \er{dV}. Then the function
$\P \in \mH_\iy(\C_+)$ and $\P $ is continuous up to the real line.
Moreover, $\p$ has a canonical factorization in $\C_+$ given by
\[
\lb{Dio} \p=\p_{in}\p_{out}.
\]
$\bu$ $\p_{in}$ is the inner factor of $\p$ having the form
\[
\lb{Di1} \p_{in}(k)={B(k)\/2ik}e^{-iK(k)},\qqq K(k)={1\/\pi}\int_\R
{d\n(t)\/k-t},\qq k\in \C_+.
\]
$\bu$ $B$ is the Blaschke product defined by \er{Bk} and $d\n(t)\ge
0$ is some singular compactly supported measure on $\R$, which
satisfies
\[
\lb{sms}
\begin{aligned}
& \n(\R)=\int_\R d\n(t)<\iy,\\
& \supp \n \ss \{z\in \R: w(z)=0\}\ss [-r_c, r_c],\qqq
r_c={\|q\|\/2}.
 \end{aligned}
\]
$\bu$  The function $K(\cdot)$ has an analytic continuation from
$\C_+$ into the domain $\C\sm [-r_c, r_c]$   and has the following Taylor series
\[
\lb{Kn} K(k)=\sum_{j=0}^\iy {K_j\/k^{j+1}},\qqq K_j={1\/\pi}\int_\R
t^jd\n(t).
\]

\no $\bu$ the function $\log |\p(t+i0)|$ belongs to $L_{loc}^1(\R)$
and $\p_{out}$ is the outer factor given by
\[
\lb{Do2} \p_{out}(k)=e^{iM(k)},\qqq  M(k)= {1\/\pi}\int_\R {\log
|\p(t)|\/k-t} dt,\qq k\in \C_+.
\]

\end{theorem}

{\bf Remark.} 1) These results are crucial to determine trace
formulas in Theorem \ref{T3}.

2) Due to  \er{asa1} the integral $M(k)$ in \er{Do2} converges
 absolutely for each $k\in \C_+$.

\medskip

We recall the well-known results. Introduce the Sobolev space $W_m$
defined  by
\[
W_m=\lt\{q\in L^1(\R): xq(x)\in L^1(\R),\  q^{(j)}\in L^1(\R),\
j=1,.., m\lt\},\qq m\geq 0.
\]
If $q\in W_{m+1}, m\ge 0$, then the function $\p(\cdot)$   satisfies
\[
\lb{apm}
  i\log \p(k)=-{Q_0\/k}-{Q_2\/k^{3}}-{Q_4\/k^{5}}+\dots
  -{Q_{2m}+o(1)\/k ^{2m+1}},
\]
as $|k |\to \iy$,  uniformly in $\arg k \in [0,\pi]$, where due to
\cite{FZ71} we have
\[
\lb{Qj} Q_0={q_0\/2}={1\/2}\int_\R q(x)dx,\qqq Q_2={1\/2^{3}}\int_\R
q^2(x)dx,...
\]
Define constants $I_j, \cJ_j$ by
\[
\begin{aligned}
\lb{Kas} I_j=\Im Q_j, \qq \cJ_0={1\/\pi}\int_0^\iy (h(t)+h(-t)) dt,\qqq
\cJ_j={1\/\pi}\int_0^\iy (h_{j-1}(t)+h_{j-1}(-t))dt,
\\
h(t)=\log |\p(t)|, \qq h_{j}=t^{j+1}(h(t)+P_{j}(t)) ,\qq
P_j(t)={I_{0}\/t}+{I_{1}\/t^2}+...+{I_{j}\/t^{j+1}}.
\end{aligned}
\]
In particular, we have $I_{2j+1}=0$ and
\[
\lb{P0}
 \cJ_1={1\/\pi}\int_0^\iy \big(t(h(t)-h(-t))+2I_0\big) dt.
 \]
The integral $\cJ_0$ converges absolutely since \er{asa1}
gives $\p(t)\p(-t)=1+{O(1)\/t^2}$ as $t\to +\iy$.

  \begin{theorem} ({\bf Trace formulas})
\lb{T3} If a potential $q$ satisfies \er{dV}, then
\[
\begin{aligned}
\lb{tr0}
& B_0+{\n(\R)\/\pi}+{1\/2}\int_\R \Re q(x)dx
={1\/\pi}\int_0^\iy \log |\p(t)\p(-t)| dt.
\end{aligned}
\]
Let a potential $q\in W_{m+1}$ for some $m\ge 0$. Then the following
identities hold true:
\[
\begin{aligned}
\lb{trj}
{B_j\/j+1}+K_j+\Re Q_j=\cJ_j,\qq j=1,...,2m,\\
\end{aligned}
\]
in particular,
\[
\begin{aligned}
\lb{tr1}
 {B_1\/2}+K_1=\cJ_1,
\end{aligned}
\]
\[
\begin{aligned}
\lb{tr2} {B_2\/3}+K_2+{1\/8}\int_\R \Re q^2(x)dx=\cJ_2.\qqq \qqq
\qqq \qqq \qqq
\end{aligned}
\]
  \end{theorem}

\no {\bf Remark.} Recall that  $B_0\ge 0$ and $K_0={\n(\R)\/\pi}\ge 0$.
Thus in order to estimate $B_0+{\n(\R)\/\pi}\ge 0$ in terms of the
potential $q$ we need to estimate the integral $\cJ_0$ in terms of
the potential $q$.

\medskip

  \begin{theorem} ({\bf Estimates})
\lb{T4} Let a potential $q$ satisfy \er{dV}. Then the following
hold true
\[
\begin{aligned}
\lb{eBs} B_0+{\n(\R)\/\pi}+{1\/2}\int_\R \Re q(x)dx \le \
{2\/\pi}\big(1+\|q\|_1\big) +{2\/\pi}r_c\big(C_0+\log r_c\big),
\end{aligned}
\]
where $\|q\|_1=\int_\R |xq(x)|dx, r_c={\|q\|\/2}$ and
$C_0=\log {3e\/2} +{3e^2\/4}+{9\/16}e^4$.
  \end{theorem}

\no {\bf Remark.} In Section~6 we  prove Theorems \ref{T3} and
\ref{T4} for the case of Schr\"odinger operators on the half-line
with the Neumann boundary condition. Recall that trace formulas of
Schr\"odinger operators on the half-line with Dirichlet boundary
condition were discussed in \cite{K18}.

Consider  estimates for  complex compactly supported potentials. In
this case the Wronskian  $w(k)$ is the entire function and has a
finite number of zeros in $\ol\C_+$.

  \begin{theorem} \lb{Tcp}
Let $q\in L^1(\R)$ and let $\supp q\ss [0,\g]$ for some $\g>0$. Then
the number of zeros $\cN(\r)$ of $w$ (counted with multiplicity) in
disk $\dD_\r(it)$ with the center $it=i2\|q\|$ and the radius $\r\ge
\sqrt 2 \|q\| $ satisfies
\[
\lb{com1} \cN(\r) \le 1+{4\/\log 2}\rt({\g \r\/\pi}+{\|q\|\/\r}\rt).
\]
In particular, the number of zeros $\cN_+$ of $w(k)$ (counted with
multiplicity) in $\ol\C_+$ satisfies
\[
\lb{com2} \cN_+ \le C_1  +C_2\g \|q\|,
\]
where the constants $C_1, C_2\le 5$ (see more about $C_1, C_2$ in
Lemma \ref{Tef}).

  \end{theorem}

Note that  the estimate of  $\cN_+$ was obtained in \cite{FLS16},
when  $q$ decays exponentially at infinity.

Our main goal is to determine trace formulas for Schr\"odinger
operators with complex potentials on the line. Our trace formula is
the identity \er{tr0},  where the left hand side is the integral
from the real part of potential, the sum  of $\Im k_j$ and the
integral of the singular measure and the right hand side is the
integral from $\log |\p(k+i0)\p(-k+i0)|$ on the real line. Here we
have  the new term, the singular measure, which is absent for real
potentials. Moreover, in \er{eBs} we estimates the singular measure
and the sum of $\Im k_j$   in terms of the potential. In our
consideration the results and technique from \cite{K17x}, \cite{K18}
are important.

In contrast to trace formulas for complex potentials, there are many
results on estimates of eigenvalues in terms of potentials, see
recent articles \cite{F18, FS17} and references therein. There exist
many recent results about bounds on sums of powers of eigenvalues
1-dimensional Schr\"odinger operators with complex-valued potentials
in terms of $L^p$-norms of the potentials published in \cite{FLLS06,
DHK09, LS09, Sa10, H11, F18}.

We shortly describe the plan of the paper. In Section~2 we present
the main properties of the Jost functions and the Wronskian $w$. In
Section~3 we prove main theorems.  Section~4 is a collection of
needed facts about Hardy spaces. In Section~5 we discuss the case of
compactly supported potentials. In Section~6 we consider
Schr\"odinger operators  on the half-line with the Neumann boundary
condition.

\section { Fundamental solutions}
\setcounter{equation}{0}

\subsection{Fundamental solutions.}
It is well known that that the Jost solution $f^+(x,k )$ of equation
\er{1.3} satisfies the  integral equation
\[
\lb{2.1} f_+(x,k)=e^{ixk}+\int_x^\iy {\sin k(t-x)\/k}q(t)f_+(t,k)dt,
\qqq  (x,k)\in [0,\iy)\ts \ol\C_+.
\]
We describe the main properties of  the Jost solution.   Due to
\er{2.1} the function $y^+(x,k)=e^{-ikx}f^+(x,k)$ satisfies the
integral equation
\[
\lb{2.2} y_+(x,k)=1+\int _x^\iy G(t-x,k)q(t)y_+(t,k)dt, \qqq
G(t,k)={\sin kt\/k}e^{ikt},
\]
$\forall\  (x,k)\in [0,\iy)\ts \ol\C_+$. The standard iterations
give $y_+(x,k)$:
\[
\lb{2.4}
\begin{aligned}
&y_+(x,k)=1+\sum_{n\ge 1}y_{+,n}(x,k),\\
& y_{+,n}(x,k)=\int _x^\iy G(t-x,k)q(t)y_{+,n-1}(t,k)dt, \qq
y_{+,0}=1.
\end{aligned}
\]
 The identity \er{2.2} gives
\[
\begin{aligned}
\lb{2.3} & f_+(0,k)=y_+(0,k)=1+\int _0^\iy {\sin kt\/k} \
q(t)f_+(t,k)dt,
\\
& f_+(0,k)'=ik-\int _0^\iy q(t)f_+(t,k)\cos kt dt .
\end{aligned}
\]
Let $q_\pm=q\c_\pm$,   where $\c_\pm$ is  the characteristic
function of the set $\R_\pm$. We recall well-known properties of the
Jost solutions  (see e.g., \cite{F63}). Define $\|q\|_1=\int_\R
|xq(x)|dx$.

\begin{lemma}
\label{TL1} Let $\int_0^\iy(1+x)|q(x)|dx<\iy$  and let $\vs_+\in
\{\|q_+\|_1,{\|q_+\|\/|k|}\}$. Then
 the functions $f_+(x,\cdot), f_+(x,\cdot)', x\ge 0$ are analytic in $\C_+$
and continuous up to the real line and satisfy
\[
\begin{aligned}
\lb{2.7}
|y_+(x,k)|\le e^{\vs_+},\\
|y_+(x,k)-1|\le \vs_+ e^{\vs_+},\\
|y_+(x,k)-1-y_1^+(x,k)|\le {\vs_+^2\/2} e^{\vs_+},
\end{aligned}
\]
and
\[
\lb{2.8}
\begin{aligned}
 |f_+(0,k)-1|\le \vs_+ e^{\vs_+},\\
|f_+(0,k)-1-f_{+,1}(0,k)|\le {\vs_+^2\/2} e^{\vs_+},\\
f_{+,1}(0,k)=\int_0^\iy {\sin kt\/k}e^{ikt}q_+(t)dt.
\end{aligned}
\]
Moreover,  $f_+(0,\cdot)'$ satisfies
\[
\lb{df1}
\begin{aligned}
 |f_+(0,k)'-ik|\le \|q_+\| e^{\vs_+},\\
|f_+(0,k)'-ik+f_{+,1}(0,k)'|\le \|q_+\|\vs_+ e^{\vs_+},\\
f_{+,1}(0,k)'=\int _0^\iy e^{ikx} \ q(x)\cos kxdx.
\end{aligned}
\]

\end{lemma}
{\bf Proof.} Let $\vs_+={\|q_+\|\/|k|}$ and
$\mD_n(x)=\{t=(t_j)_1^n\in \R^n: x=t_0<t_1< t_2<...< t_n\}$ for
$x>0$. Substituting the estimate $ |G(t,k)|\le {1\/|k|}$ for all
$t\ge 0, \ k\in \ol\C_+\sm\{0\}$ into the identity
\[
\lb{ynx} y_{+,n}(x,k)=\int\limits_{\mD_n(x)} \lt(\prod\limits_{1\le
j\le n} G(t_{j}-t_{j-1},k)q(t_j)dt_j\rt),
\]
we obtain
\[
\begin{aligned}
\lb{2.10} |y_{+,n}(x,k)|\le {1\/|k|^n}\int\limits_{\mD_n(x)}
\lt(\prod\limits_{1\le j\le n} |q(t_j)|dt_j\rt) = {\|q_+\|^n\/n!
|k|^n}.
\end{aligned}
\]
This shows that for each $x\ge 0$ the series \er{2.4} converges
uniformly on bounded subset of $\ol\C_+\sm \{|k|>\ve\}$ for any
$\ve>0$. Each term of this series is an analytic function in $\C_+$.
Hence the sum is an analytic function in $\C_+$. Summing the
majorants we obtain estimates  \er{2.7}-\er{2.8} for
$\vs_+={\|q_+\|\/|k|}$. Thus the functions $f_+(x,\cdot),
{f_+}'(x,\cdot),x\ge 0$ are analytic in $\C_+$ and continuous up to
the real line without the point $0$.

Let in addition $\|q_+\|_1=\int_0^\iy x|q(x)|dx<\iy$ and let
$\vs_+=\|q_+\|_1$. The function $G(t,k)={\sin kt\/k}e^{ikt}$ satisfy
$|G(t,k)|\le t$ for all $k\in \ol \C_+, t>0$. Then using above
arguments
 we obtain
\[
\begin{aligned}
\lb{2.15} |y_{+,n}(x,k)|\le \int\limits_{\mD_n(x)}
\lt(\prod\limits_{1\le j\le n} |(t_{j}-t_{j-1})q(t_j)|dt_j\rt)
\\
\le \int\limits_{\mD_n(x)} \lt(\prod\limits_{1\le j\le n}
|t_{j}q(t_j)|dt_j\rt)={\|q_+\|_1^n\/n! }.
\end{aligned}
\]
This shows that for each $x\ge 0$ the series \er{2.4} converges
uniformly on bounded subset of $\ol\C_+$. Each term of this series
is an analytic function in $\C_+$. Hence the sum is an analytic
function in $\C_+$. Summing the majorants we obtain estimates
\er{2.7}-\er{2.8} for $\vs_+=\|q_+\|_1$. Thus the functions
$f_+(x,\cdot), f_+(x,\cdot)'\, x\ge 0$ are analytic in $\C_+$ and
continuous up to the real line.

Consider $f_+(0,k)'$. From \er{2.7} and \er{2.3}  we get the first
estimate in \er{df1}.  From \er{2.3} we have
\[
\begin{aligned}
\lb{fx}
 f_+'(0,k)-ik=-\int _0^\iy (\cos kx)e^{ikx} \ q(x)y_+(x,k)dx
 =-f_{+,1}(0,k)'-F_2(k),\\
 f_{+,1}'(0,k)=\int _0^\iy (\cos kx )e^{ikx} \ q(x)dx,\qq
 F_2(k)=\int _0^\iy (\cos kx) e^{ikx} q(x)(y_+(x,k)-1)dx,
\end{aligned}
\]
where
$$
|F_2(k)|\le \int _0^\iy |q(x)| |y_+(x,k)-1|dx\le
\|q_+\|\vs_+e^{\vs_+},
$$
which yields \er{df1}. \BBox

We consider the Jost function $f_-(x,\cdot)$, defined by \er{pf},
which satisfies the integral equation
\[
\lb{fpx} f_-(x,k)=e^{-ikx}-\int_{-\iy}^x {\sin
k(t-x)\/k}q(t)f_-(t,k)dt, \qqq  x\le 0,\qq k\in \ol\C_+,
\]
and the function $y_-(x,k)=e^{ikx}f_-(x,k)$ also satisfies the
integral equation
\[
\lb{24p}
\begin{aligned}
&y_-(x,k)=1-\int_{-\iy}^x {\sin k(t-x)\/k}q(t)e^{ik(x-t)}y_-(t,k)dt.
\end{aligned}
\]

\begin{lemma}
\label{TL2} Let $\int_{-\iy}^0 (1+|x|)|q(x)|dx<\iy$ and let $x\le 0$
 and let $\vs_-\in \rt\{\|q_-\|_1,{\|q_-\|\/|k|}\rt\}$. Then
the functions $f_-(x,\cdot), f_-(x,\cdot)'$  are analytic in $\C_+$
and continuous up to the real line and satisfy
\[
\begin{aligned}
\lb{2.7x}
|y_-(x,k)|\le e^{\vs_-},\\
|y_-(x,k)-1|\le \vs_- e^{\vs_-},\\
|y_-(x,k)-1-y_1^-(x,k)|\le {\vs_-^2\/2} e^{\vs_-},
\end{aligned}
\]
and
\[
\lb{2.8x}
\begin{aligned}
 |f_-(0,k)-1|\le \vs_- e^{\vs_-},\\
|f_-(0,k)-1-f_{1}^-(0,k)|\le {\vs_-^2\/2} e^{\vs_-},\\
f_{-,1}(0,k)=-\int_{-\iy}^0 {\sin kt\/k}e^{-ikt}q_-(t)dt,
\end{aligned}
\]
and ${f_-}'(0,\cdot)$ satisfies
\[
\lb{2.8xy}
\begin{aligned}
 |f_-(0,k)'-ik|\le \|q_-\| e^{\vs_-},\\
|f_-(0,k)'-ik+f_{-1}'(0,k)|\le \|q_-\| \vs_- e^{\vs_-},\\
f_{-,1}'(0,k)=\int_{-\iy}^0 q(x)e^{ikx}\cos kx \ dx
\end{aligned}
\]
\end{lemma}
{\bf Proof.} The proof repeats the case of Lemma \ref{TL1}. \BBox

Using Lemmas \ref{TL1} and  \ref{TL2} we have the following decompositions:
\[
\lb{f12}
\begin{aligned}
 f_+(0,k)=1+F_s^+(k),\qqq F_s^+(k)=\int_0^\iy {\sin kt\/k}q(t)f_+(t,k)dt, \\
f_+(0,k)'=ik-kF_c^+(k),\qqq F_c^+(k)=\int_0^\iy {\cos kt\/k}q(t)f_+(t,k)dt, \\
\end{aligned}
\]
and
\[
\lb{f123}
\begin{aligned}
 f_-(0,k)=1-F_s^-(k),\qqq F_s^-(k)=\int_{-\iy}^0 {\sin kt\/k}q(t)f_-(t,k)dt, \\
{f_-}'(0,k)=-ik+kF_c^-(k),\qqq F_c^-(k)=\int_0^\iy {\cos
kt\/k}q(t)f_-(t,k)dt.
\end{aligned}
\]

\begin{lemma}
\label{TL3} Let $\int_\R (1+|x|)|q(x)|dx<\iy$ and let $\vs\in
\rt\{\|q\|_1,{\|q\|\/|k|}\rt\}$. Then the Wronskian $w$ is analytic
in $\C_+$ and continuous up to the real line and satisfies
\[
\begin{aligned}
\lb{w1}
|w(k)|\le (2|k|+\|q\|)e^{\vs},\\
\end{aligned}
\]
and
\[
\lb{w2}
\begin{aligned}
& w=2ik-w_1(k)-w_2(k),\\
& w_1(k)=\int_0^{\iy}q(t)y_+(t,k)dt+\int_{-\iy}^0 q(t)y_-(t,k)dt,\\
& w_2(k)= k(F_s^+F_c^--F_c^+F_s^-)=\int_0^\iy dx
\int_{-\iy}^0{(1-e^{i2k(x-y)})\/2ik}q(x)q(y)y_+(x,k)y_-(y,k)dy,
\end{aligned}
\]
and
\[
\begin{aligned}
\lb{w3} & |w_1(k)|\le \|q\|e^{\vs},\qqq  |w_1(k)-q_0|\le
\|q\|e^{\vs},
 \\
&  |w_2(k)|\le (\|xq_+\|\|q_-\|+\|q_+\|\|xq_-\|) e^{\vs},\qqq
|w_2(k)|\le {\|q_-\|\|q_+\|\/|k|}e^{\vs},
\end{aligned}
\]
where $q_0=\int_\R qdx$, and  we have at $\vs={\|q\|\/|k|}$:
\[
\lb{w4}
\begin{aligned}
|\p(k)-1|\le {(\vs+\vs^2)\/2} e^{\vs},
\\
 |\p(k)-1+{q_0\/2ik}|\le  \vs^2 e^{\vs}.
\end{aligned}
\]

\end{lemma}
{\bf Proof.} We show \er{w1}. From  Lemmas \ref{TL1}, \ref{TL2} we have
$$
\begin{aligned}
 |w(k)|=|f_-(0,k)f_+(0,k)'-f_-(0,k)'f_+(0,k)|\\
 \le
 e^{\vs_-}(|k|+\|q_+\| e^{\vs_+})  +(|k|+\|q_-\| e^{\vs_-})e^{\vs_+}
\le (2|k|+\|q\|)e^{\vs}.
\end{aligned}
$$
We show \er{w2}. Let $q\in L^1(\R)$. Using \er{f12}, \er{f123} we
obtain
\[
\lb{f13}
\begin{aligned}
w=\{f_-, f_+\}=(f_-{f_+}'-{f_-}'f_+)(0,k)
\\
=(1-F_s^-)(ik-kF_c^+)-(-ik+kF_c^-)(1+F_s^+)
\\
=2ik-k(F_c^+-iF_s^+)-k(F_c^-+iF_s^-)-k(F_s^+F_c^--F_c^+F_s^-),
\end{aligned}
\]
where
\[
\lb{f14}
\begin{aligned}
k(F_c^+-iF_s^+)=\int_0^{\iy}e^{-ikt}q(t)f_+(t,k)dt=\int_0^{\iy}q(t)y_+(t,k)dt,
 \qqq  \\
k(F_c^-+iF_s^-)=\int_0^\iy e^{ikt}q(t)f_-(t,k)dt=\int_0^\iy
q(t)y_-(t,k)dt.
\end{aligned}
\]

Thus we have
\[
\lb{f15}
\begin{aligned}
w=2ik-\int_0^{\iy}q(t)y_+(t,k)dt-\int_{-\iy}^0
q(t)y_-(t,k)dt-k(F_s^+F_c^--F_c^+F_s^-).
\end{aligned}
\]

Let $c_x={\cos kx}$ and $ s_x={\sin kx}$. We get
\[
\lb{f16}
\begin{aligned}
w_2=k(F_s^+F_c^--F_c^+F_s^-)=
\int_0^\iy dx \int_{-\iy}^0{(s_xc_y-c_xs_y)\/k}q(x)q(y)f_+(x,k)f_-(y,k)dy\\
=-\int_0^\iy dx
\int_{-\iy}^0{(1-e^{i2k(x-y)})\/2ik}q(x)q(y)y_+(x,k)y_-(y,k)dy,
\end{aligned}
\]
which yields \er{w2}.
We show \er{w3}. Let $q_0=\int_\R q(x)dx$.  From  \er{2.7}, \er{2.7x} we have
$$
\begin{aligned}
w_1(k)=q_0+\int_0^{\iy}q(x)(y_+(x,k)-1)dx+\int_{-\iy}^0
q(x)(y_-(x,k)-1)dx,\\
|w_1(k)-q_0|\le\int_0^{\iy}|q(x)(y_+(x,k)-1)|dt+\int_{-\iy}^0
|q(x)(y_-(x,k)-1)|dx\\
\le\|q_+\|\vs_+e^{\vs_+}+\|q_-\|\vs_-e^{\vs_-}\le \|q\|\vs e^{\vs}.
\end{aligned}
$$
Similar arguments yield  $|w_1(k)|\le \|q\|e^{\vs}$.

From  \er{2.7}, \er{2.7x} if $|k|\ge 1$ we obtain
$$
\begin{aligned}
|w_2(k)|\le \int_0^\iy dx
\int_{-\iy}^0|q(x)q(y)y_+(x,k)y_-(y,k)|{dy\/|k|}\le
{\|q_+\|\|q_-\|\/|k|} e^{\vs_++\vs_-}\le \vs \|q\|e^{\vs},
\end{aligned}
$$
and if $|k|\le 1$  we get
$$
\begin{aligned}
|w_2(k)|\le \int_{\R_+}dx
\int_{\R_-}(x-y)|q(x)q(y)y_+(x,k)y_-(y,k)|dy \le
(\|xq_+\|\|q_-\|+\|q_+\|\|xq_-\|) e^{\vs},
\end{aligned}
$$
where the simple estimate has been used :
$|1-e^{ikz}|\le z|k|$  for all  $(k,z)\in \C_+\ts \R_+$.

We show \er{w4}. From \er{w3} we have $w(k)-2ik+q_0=(q_0-w_1)-w_2$
and
$$
\begin{aligned}
 |w(k)-2ik+q_0|=|(q_0-w_1)-w_2| \le 2\|q\|\vs e^{\vs},
\end{aligned}
$$
which yields the first estimate in \er{w4}. Similar arguments give
the second one in \er{w4}. \BBox

\section { Proof of main theorems }
\setcounter{equation}{0}


In order to study zeros of the function $\p(k)={w(k)\/2ik}$ in the
upper-half plane we need to study the Blaschke product, defined by
\er{Bk}. Recall that in order to describe the basic properties of
the Blaschke product $B$ as an analytic function in $\C_+$ we modify
the function $\p$ and define the modified function by
$\P(k)={w(k)\/2i(k+i)},\ k\in \C_+$.
 We recall the well-known identity
\[
\lb{dD}
 \p(k)=\det (I+Y_0(k)), \qqq k\in \C_+,
\]
where that $Y_0(k)$ is a trace class operator given by
$$
Y_0(k)=|q|^{1\/2} R_0(k)|q|^{1\/2}e^{i\arg q},\qqq
R_0(k)=(H_0-k^2)^{-1},\qqq k\in \C_+.
$$

{\bf  Proof Proposition \ref{T1}.}  The zeros of $w, \p$ and $\P$ in
$\C_+$ are the same. Due to estimate \er{w1} the function $\P\in
\mH_\iy$.
The estimate \er{w4} gives that all zeros of $\P$ are uniformly
bounded. Note that (see page 53 in \cite{G81}), in general, in the upper half
plane the condition \er{B1} is replace by
\[
\lb{BLy} \sum {\Im k_j\/1+|k_j|^2}<\iy,
\]
and the Blaschke product with zeros $k_j$ has the form
\[
\lb{BL2x} B(z)={(k-i)^m\/(k+i)^m} \prod_{k_j\ne
0}^N{|1+k_j^2|\/1+k_j^2}\rt(\frac{k-k_j}{k-\ol k_j}\rt), \qqq k\in
\C_+.
\]
If all  moduli $|k_n|$ are uniformly bounded,  the estimate \er{BLy}
becomes $\sum \Im k_j<\iy$ and the convergence factors in \er{BL2x}
are not needed, since $\prod_{k_j\ne 0}^N\big(\frac{k-k_j}{k-\ol
k_j}\big)$ already converges.

The statement i) is a standard fact for the function $\P\in
\mH_\iy$, see Sect. VI  in \cite{Ko98}. Lemma \ref{TL3} and Lemma
\ref{Tf1} imply ii). \BBox

We describe the determinant $\p(k), k\in \C_+$ in terms of a
canonical factorization.

{\bf Proof of Theorem \ref{T2}.} From Proposition \ref{T1} we have that
the modified function $\P\in\mH_\iy$
and due to \er{w4}the function $\P$ has asymptotics
$\P(k)=1+O(1/k)$ as $|k|\to \iy$ uniformly in $\arg k \in [0,\pi]$.
Then from Theorem \ref{Tcf} we deduce that  the function $\P\in \mH_\iy$
has a canonical factorization in $\C_+$ given by
\[
\begin{aligned}
\lb{P3} \P=\P_{in}\P_{out},
\qqq
\P_{in}(k)=B(k)e^{-iK(k)},\qqq  K(k)={1\/\pi}\int_\R {d\n(t)\/k-t}.
\end{aligned}
\]
$\bu$  $d\n(t)\ge 0$ is some singular compactly supported measure on
$\R$, which satisfies
\[
\lb{P4}
\begin{aligned}
\n(\R)=\int_\R d\n(t)<\iy,\qqq
 \supp \n \ss \{k\in \R: \P(k)=0\}\ss [-r_c, r_c],
 \end{aligned}
\]
for $r_c={\|q\|\/2}$, since $w$ satisfies \er{1.2}.

$\bu$  The function $K(\cdot)$ has an analytic continuation from
$\C_+$ into the domain $\C\sm [-r_c, r_c]$   and has the following
Taylor series
\[
\lb{P5} K(k)=\sum_{j=0}^\iy {K_j\/k^{j+1}},\qqq K_j={1\/\pi}\int_\R
t^jd\n(t).
\]

$\bu$ $B$ is the Blaschke product for $\Im k>0$ given by \er{Bk}.

 $\bu$ $\P_{out}$ is the outer factor given by
 $\P_{out}(k)=e^{iM(k)}$, where
 $ M(k)= {1\/\pi}\int_\R {\log |\P(t)|\/k-t} dt, k\in \C_+$.

 We consider the function $\P=\p \x$, where $\x={k\/k+i}$.
  It is clear that $\x={k\/k+i}$ has the following factorization
\[
\lb{P2} \x(k)={k\/k+i}=\exp {{1\/\pi}\int_\R { \log |\x(t)|dt\/k-t} },
\qqq k\in \C_+,
\]
where $ \log \x(k)$ in $\C_+$ is defined by $\log \x(k)=O(1/k)$ as $|k|\to \iy$.
This yields
$$
M(k)= {1\/\pi}\int_\R {\log |\p(t)|\/k-t} dt+{1\/\pi}\int_\R {\log
|\x(t)|\/k-t} dt={1\/\pi}\int_\R {\log |\p(t)|\/k-t} dt+\log \x(k),
$$
for all $k\in \C_+$.
Here $\log |\P|, \log |\x|\in L_{loc}^1(\R)$, which yields $\log
|\p|\in L_{loc}^1(\R)$, since $\P=\p \x$.
Thus from the properties of $\P$ we  obtain all properties of $\p$
 formulated in Theorem
\ref{T2}. \BBox

Remark that  a canonical factorization is the first trace formula.
It is a generating function,  the differentiation of a canonical
factorization produces an identity for $\Tr (R(k)-R_0(k)), k\in
\C_+$, where $R(k)=(H-k^2)^{-1}$ and $R_0(k)=(-\pa_x^2-k^2)^{-1}$ is
the free resolvent.

\begin{corollary}
\lb{T2x} Let a potential $q$ satisfy \er{dV}. Then  the trace
formula
\[
\lb{tre1} -2k\Tr \rt(R(k)-R_0(k) \rt)= \sum {2i\Im
k_j\/(k-k_j)(k-\ol k_j)}+{i\/\pi}\int_{\R}{d\m(t)\/(t-k)^2},
\]
holds true for any $k\in \C_+\sm\dk_q$, where  the measure
$d\m(t)=\log |\p(t)|dt-d\n(t)$ and the series converges uniformly in
every bounded disc in $\C_+\sm \dk_q$.

\end{corollary}

{\bf Proof}. We repeat arguments from the proof of Corollary 1.3
from \cite{K17}. Differentiating \er{Dio} and using Theorem \ref{T2}
we obtain
\[
\begin{aligned}
\lb{derfx}
{\p'(k)\/\p(k)}={B'(k)\/B(k)}-{i\/\pi}\int_\R{d\m(t)\/(k-t)^2},\qq
\qqq {B'(k)\/B(k)}=\sum {2i\Im k_j\/(k-k_j)(k-\ol k_j)},\qq \forall
\ k\in \C_+,
\end{aligned}
\]
where $d\m(t)=h(t)dt-d\n(t)$. The derivative of the determinant
$\p=\det (I+Y_0(k))$ defined by \er{dD} satisfies
\[
\lb{derfz} {\p'(k)\/\p(k)}=-2k\Tr \big(R(k)-R_0(k)\big), \qq \forall
\ k\in \C_+,
\]
see \cite{GK69}. Combining \er{derfx}, \er{derfz} we obtain
\er{tre1}. Note that the series converges uniformly in every bounded
disc in $\C_+\sm \dk_q$, since $\sum \Im k_j<\iy$. \BBox

We prove the first main result about the trace formulas.

\no {\bf  Proof of Theorem \ref{T3}.} Let a potential $q$ satisfy
\er{dV}. Then due to Lemma \ref{TL1} the function  $\p(t)=1-{\Im
Q_0+O(1/t)\/t}$ as $t\to \pm \iy$ and due to Theorem \ref{Tcf} the
function $h(t)=\log|\p(t)|, t\in \R$, belongs to
 $L^1_{loc}(\R)$. Then the function $h(t)$ satisfies all conditions in
Lemma \ref{Th1} for $m=0$ and the function $M(k)={1\/\pi}\int_\R
{h(t)\/k-t}dt,k\in \C_+$ satisfies $M(k)={\cJ_0+iI_{0}+o(1)\/k}$ as
$k=iv, v\to \iy$. From Lemma \ref{Th1} and from asymptotics \er{apm}
we obtain
$$
\begin{aligned}
\exp \rt[i{Q_0+o(1)\/k }  \rt] =\exp \rt[ -i{B_0\/k}-i{K_0\/k}
+i{\cJ_0+iI_{0}\/k} +{o(1)\/k}\rt],
\end{aligned}
$$
as $ k=iv, \ v\to \iy$, which yields $ \Re Q_0=B_0+K_0-\cJ_0$, since
$I_j=\Im Q_j$. Thus we have \er{tr0}. The proof of \er{trj} is
similar. \BBox

\no {\bf  Proof of Theorem \ref{T4}.} We estimate the integral
$\cJ_0={1\/\pi}\int_0^\iy \x(k) dk$ in \er{P0}, where $\x(k)=\log |\p(k)\p(-k)|$.
 We rewrite $\cJ_0$ in the following form
\[
\lb{F1}
\begin{aligned}
\cJ_0
={1\/\pi}(\cJ_{01}+\cJ_{02}),\qqq \cJ_{01}=\int_0^{\ve} \x(k) dk,\qqq
\cJ_{02}=\int_{\ve}^\iy \x(k) dk,
 \end{aligned}
 \]
where $\ve:=\|q\|$ for shortness. Consider $\cJ_{01}$ for  the case
$\ve\ge1$. We have
\[
\lb{J01x} \cJ_{01}=\int_0^1\x(k)dk+\int_1^{\ve}\x(k) dk.
\]
The estimate \er{w1} gives with $\vs={\|q\|_1}$ for $k\in (0,1)$:
\[
\lb{J01}
\begin{aligned}
\int_0^1\x(k)dk\le 2\int_0^1 \log \rt({k+r_c\/k}e^\o\rt)
dk={2}\|q\|_1+{2}\int_0^1\log
\rt(1+{r_c\/k}\rt)dk\\
={2}\|q\|_1+{2}\log(1+r_c)+{2}\int_0^1{{r_c\/k}dk\/1+{r_c\/k}}
\le {2}\|q\|_1+{2}\log(1+r_c)+{2},
 \end{aligned}
\]
and with $\vs={\|q\|\/|k|}$ for $k\in (1,\ve)$:
$$
\begin{aligned}
 \int_1^{\ve}\x(k) dk\le {2\|q\|}\int_1^{\ve}
{dk\/k}+{2}\int_1^{\ve}\log \rt(1+{r_c\/k}\rt)dk
\\
=2\ve \log \ve +2k\log
\rt(1+{r_c\/k}\rt)\rt|_1^{\ve }+{2}\int_1^{\ve }{{r_c\/k}dk\/1+{r_c\/k}}
\\
\le 2\ve \log \ve+{2}\ve\log {3\/2}-{2}\log (1+r_c)+2(\ve-1).
 \end{aligned}
$$
This estimate and \er{J01} give
\[
\lb{J01z}
\cJ_{01}\le {2}\|q\|_1+{2}\ve\big(1+\log {3\/2}+\log \ve\big).
\]
Consider $\cJ_{02}$. Due to \er{w3}  we define $g(k), \wt g(k)$ by
$\p=1+g$ and $g=-{q_0\/2ik}+\wt g(k)$. The estimate \er{w3} gives
for $k>\ve$:
\[
\begin{aligned}
\lb{F1x}
 |\wt
g(k)|\le \vs^2 e^{\vs},\qqq |\p(k)-1|\le {1\/2}\vs(1+\vs) e^{\vs}\le
\vs e^{\vs},\\
\end{aligned}
 \]
where $\vs={\|q\|\/|k|}<1$. This yields
\[
\lb{asxd}
\begin{aligned}
f(t)=\p(t)\p(-t)=(1+g(k))(1+g(-k))=1+g(k)+g(-k)+g(k)g(-k),\\
\end{aligned}
\]
where
\[
\lb{asp1Vxx}
\begin{aligned}
& |g(k)g(-k)|\le \vs^2 e^{2\vs},
\\
& g(k)+g(-k)=\wt g(k)+\wt g(-k), \qqq |\wt g(k)+\wt g(-k)|\le
2{\vs^2} e^{2\vs}.
\end{aligned}
\]
Thus we obtain
\[
\lb{ah1}
\begin{aligned}
f(k)=1+f_1(k),\qqq f_1(k)=\wt g(k)+\wt g(-k)+g(k)g(-k),\qq
 |f_1(k)|\le 3{\vs^2} e^{2\vs}.
\end{aligned}
\]
This yields
\[
\lb{ah1x}
\begin{aligned}
f\ol f=(1+f_1)(1+\ol f_1)=1+2\Re f_1+|f_1|^2,
\end{aligned}
\]
and then
\[
\begin{aligned}
\log |f|^2\le F, \qq F=2 \Re f_1+|f_1|^2\le 6{\vs^2} e^{2\vs}+9\vs^4
e^{4\vs}.
\end{aligned}
\]

We have the identities
$$
\begin{aligned}
& V_m=\int_{\ve}^\iy {\ve^m\/t^m} e^{m\ve \/t}dt=
\ve \int_{1}^\iy {1\/t^m} e^{m\/t}dt=
\ve \int_0^{1} s^{m-2} e^{ms}ds,
\\
& V_2=\ve {e^2-1\/2},\qqq V_4=\ve{5e^4-1\/32}.
\end{aligned}
$$
Then we obtain
\[
\begin{aligned}
\cJ_{02}={1\/2\pi}\int_{\ve}^\iy \log |f(k)|^2 dk\le
{1\/2\pi}\int_{\ve}^\iy Fdk,
\end{aligned}
\]
and
\[
\lb{JF}
\begin{aligned}
\int_{\ve}^\iy Fdk\le \int_{\ve}^\iy(6{\vs^2} e^{2\vs}+9\vs^2
e^{4\vs})dk=6V_2+9V_4 =\ve C_1,
\end{aligned}
\]
where $C_1=3(e^2-1)+{9\/32}(e^4-1)$. Thus collecting
\er{J01z}-\er{JF} we obtain
\[
\lb{F1xx}
\begin{aligned}
\cJ_0\le  {2\|q\|_1\/\pi}+{2\ve\/\pi}
\big(1+\log {3\/2}+\log \ve\big)+{\ve \/\pi}C_1
={{2}\|q\|_1\/\pi}+{2\ve\/\pi}\big(1+C_1+\log {3\/2}+\log \ve\big).
 \end{aligned}
 \]

Consider the case $\ve<1$. We need to estimate only $\cJ_{01}$,
since we obtain the estimate $\cJ_{02}$ for any $\ve >0$. The
estimate \er{w1} gives with $\vs={\|q\|_1}$ for $k\in (0,\ve)$:
\[
\lb{J01a}
\begin{aligned}
\int_0^\ve\x(k)dk\le 2\int_0^\ve \log \rt({k+r_c\/k}e^\o\rt)
dk={2}\|q\|_1\ve+{2}\int_0^\ve\log \rt(1+{\ve\/2k}\rt)dk
\\
={2}\ve\|q\|_1+{2}\ve \int_0^1\log
\rt(1+{1\/2t}\rt)dt={2}\ve\|q\|_1+{2}\ve \log{3^{3\/2}\/2},
 \end{aligned}
\]
since $\int_0^1\log\big(1+{1\/2t}\big)dt=\log{3^{3\/2}\/2}$.
Thus collecting \er{J01a}, \er{JF} we obtain
$$
\begin{aligned}
\cJ_0\le  {2\ve\/\pi}\|q\|_1+{2\ve\/\pi}
\log{3^{3\/2}\/2}+{\ve\/\pi} C_1.
 \end{aligned}
$$

\BBox




\

\section { Appendix: Analytic functions in the upper half-plane}
\setcounter{equation}{0}

\

Let $f\in \mH_\iy(\C_+)$ and let $\{k_j\}$ be its zeros, uniformly
bounded by $r_0$. Then its  Blaschke product $B\in \mH_\iy(\C_+)$
and satisfies
\[
\lb{B3} \lim_{v\to+0} B(u+iv)=B(u+i0), \qqq |B(u+i0)|=1 \qqq   \
{\rm almost\ e.w.\ for} \ u\in\R,
\]
\[
\lb{B4} \lim _{v\to 0}\int_\R\log |B(u+iv)|du=0.
\]
The function $\log B(k)$ is analytic in $\{|k|>r_0\}$ and has the
corresponding Tailor series given by
\[
\lb{Ba3} \log
B(k)=-{iB_0\/k}-{iB_1\/2k^2}-{iB_2\/3k^3}-....-{iB_{n-1}\/nk^n}-....
\]
where each sum $B_n=2\sum_{j} \Im k_j^{n+1}, \ n\ge 0$ is bounded
and satisfies
\[
\lb{Ba1} |B_n|\le 2\sum |\Im k_j^{n+1}|\le {\pi}(n+1)r_0^{n}B_0<\iy
\qqq \qqq \forall \ n\ge 1,
\]
see e.g., \cite{K17x}. We recall the standard facts about the
canonical factorization, see e.g. \cite{G81}, \cite{Ko98} and  in
the needed specific  form for us from \cite{K17x}.

\begin{theorem}\label{Tcf}
Let a function $f\in\mH_p$  and $f(t+i0)=1+O(t^{-a})$ as $t\to\pm
\iy$ for some $p\ge 1, a>0$.  Then $f$ has a canonical factorization
in $\C_+$ given by
\[
\begin{aligned}
\lb{af1} f=f_{in}f_{out},\qqq f_{in}(k)=B(k)e^{-iK(k)},\qqq
K(k)={1\/\pi}\int_\R {d\n(t)\/k-t}.
\end{aligned}
\]
$\bu$  $d\n(t)\ge 0$ is some singular compactly supported measure on
$\R$, which satisfies
\[
\lb{smsx}
\begin{aligned}
\n(\R)=\int_\R d\n(t)<\iy,\\
 \supp \n \ss \{k\in \R: f(k)=0\}\ss [-r_c, r_c],
 \end{aligned}
\]
for some $r_c>0$.

$\bu$  The function $K(\cdot)$ has an analytic continuation from
$\C_+$ into the domain $\C\sm [-r_c, r_c]$   and has the following Taylor series
\[
\lb{Knx} K(k)=\sum_{j=0}^\iy {K_j\/k^{j+1}},\qqq K_j={1\/\pi}\int_\R
t^jd\n(t).
\]

$\bu$ $B$ is the Blaschke product for $\Im k>0$ given by \er{Bk}.

 $\bu$ $f_{out}$ is the outer factor given by
\[
\lb{Do2x} f_{out}(k)=e^{iM(k)},\qqq  M(k)= {1\/\pi}\int_\R {\log
|f(t)|\/k-t} dt,\qq k\in \C_+,
\]
where the function $\log |f(t+i0)|$ belongs to  $L_{loc}^1(\R)$.
\end{theorem}

{\bf Remark.} 1) These results are crucial to determine trace
formulas in Theorem \ref{T3}.

2) The integral $M(k)$ in \er{Do2x} converges absolutely, since
$f(t)=1+{O(1)\/t^a}$ as $t\to \pm\iy$.

\medskip

In order to determine trace formulas in Theorem \ref{T3} we recall
asymptotics from \cite{K17x}.

\begin{lemma}
\lb{Th1}

Let $h$ be a real function from  $L_{loc}^1(\R) $ and let
 $h$ satisfy
\[
\lb{hm} h(t)=-P_m(t)-{O(1)\/t^{m+2}},\qq
P_m(t)={I_{0}\/t}+{I_{1}\/t^2}+...+{I_{m}\/t^{m+a}},
\]
as \ $t\to \pm \iy$, for some constants $a>0$, $I_0,I_1,....,I_m\in
\R$ and integer $m\ge 0$. Then
\[
\begin{aligned}
\lb{asM}
{1\/\pi}\int_\R {
h(t)\/k-t}dt={J_0+iI_{0}\/k}+{J_1+iI_{1}\/k^2}+...+{J_m+iI_{m}+o(1)\/k^{m+1}},\\
\end{aligned}
\]
at $k=iv, v\to\iy$, where
$$
\begin{aligned}
J_j={1\/\pi}\int_0^\iy (h_j(t)+h_j(-t))dt,\qq h_0=h,\qq
 h_{j}=t^j(h(t)+P_{j-1}(t)),\qquad j=0,1,2,...,m.
\end{aligned}
$$
\end{lemma}

\begin{lemma}
\lb{Tf1} Let  a function $f\in \mH_\iy$ and be continuous up to the
real line. Assume that it satisfies
\[
\lb{fx1}
\begin{aligned}
|f(k)-(2ik-\t_o)|\le \t \vs e^{\vs},\qqq \forall \ k\in \ol\C_+, \qq
\vs= \min \{\s, {\t\/|k|}\},
\end{aligned}
\]
for some positive constants $\s,\t$ and $\t_o\in \C \sm \R$ and
$|\t_o|\le \t$. Then

 i) If $\t \s e^{\s}<\Re \t_o$, then the function $f$
does  not have zeros in $\ol\C_+$.

ii)  If $\t \s e^{\s}<-\Re \t_o$, then function $f$ has exactly  one
simple zeros in $\ol\C_+$.

\end{lemma}
{\bf Proof.}  Define the function $\wt f(k)=f(k)-(2ik -\t_o)$, where
$|\wt f(k)|\le \t \s e^{\s}.$ We define the new variable
$t={2k\/\t}\in \ol\C_+$ and rewrite $f$ in the form $f(k)=i\t g(t)$,
where
\[
\lb{f2}
\begin{aligned}
g(t)=t+i\ve+i\wt g(t),\qq \qq \ve={\t_o\/\t},\qq \wt g(t)={\wt
f(k)\/\t}, \qq |\wt g(t)|\le \s e^{\s},\qqq |\ve|\le 1.
\end{aligned}
\]

i)  Consider the first case, and let $\s e^{\s}<\Re \ve$. Then we
have for all $t\in \ol\C_+$:
\[
\lb{f3} \Im g(t)=\Im t+\Re \ve+\Re\wt g(t)\ge \Im t+\Re \ve-\s
e^{\s}>0.
\]

ii) Consider the second case, let $\ve=-r +i\n, r>0$ and let $\s
e^{\s}<-\Re \ve=r$. Then $i\ve =-(\n+ir)$ and define the new
variable $z=t-\n$ and
$$
g=z-i r  +\wt g(t), \qqq g_0=z-i r.
$$
Introduce the contour $\cC_\r=\cC_\r^+\cup I_\r$, where $ \cC_\r^+=
\{|z|=\r, z\in \ol\C_+\}$ and $I_\r=[-\r,\r]$ for $\r\ge 2r$.
 Then,  we have for all $z\in \cC_\r$ and any $\r\ge 2r$:
$$
|g(t)-g_0(z)|=|\wt g(t)| \le \s e^{\s}=|g_0(z)|
{ue^{u}\/|g_0(z)|}<|g_0(z)|,  \qq t=z+\n,
$$
where we have used ${\s e^{\s}\/|g_0(z)|}\le {\s e^{\s}\/r}<1$ on
$\cC_\r$, since $ |g_0(t)|\ge r$ for all \ $t\in \R, $ and
$$
|g_0(t)|=|\r e^{i\f}-ir|\ge \r-r\ge r\qq \forall \ t\in \C_\r.
$$
Thus, by Rouche's theorem, $g$ has one simple root,  as $g_0$ in the
 region $\{|z|\le\r, z\in \ol\C_+\}$.
\BBox


\section {Schr\"odinger operators with compactly supported potentials}
\setcounter{equation}{0}

\subsection{Entire functions}
An entire function $f(k)$ is said to be of exponential type if there
is a constant $\b$ such that $|f(k)|\leq\const e^{\b |k|}$
everywhere. The infimum of the set of $\b$ for which such inequality
holds is called the type of $f$.

\no {\bf Definition.} {\it Let $\cE_\g, \g>0$ denote the class of
exponential type  functions $f$ satisfying
\[
\lb{dE}
\begin{aligned}
& |f(k)-(2ik-\t_o)|\le \t \vs e^{2\g k_-+\vs},\qqq \forall \ k\in
\ol\C_+, \qq \vs(k)= \min \{\s, {\t\/|k|}\},
\end{aligned}
\]
where $k_-={1\/2}(|\Im k|-\Im k)\ge 0$,  for some constants $(\s,\t,
\t_o)\in \R_+^2\ts \C$ and $|\t_o|\le v$ and
 each its zero $z$ in $\ol\C_+$ satisfies $|z|\le {\t\/2}$.
}

Note that if $q\in L^1(\R_+)$ and $\supp q\ss [0,\g]$, then the
Wronskian $w\in \cE_\g$, see below the proof of Theorem \ref{Tcp}.
Define the disk $\dD_r(t)=\{z: |t-z|<r\}$ for $t\in \C$ and $r>0$.

\begin{lemma}
\label{Tef} i)  Let $f\in \cE_\g, \g>0$.    If $k\in \C_+ $ and
$|k|\ge 2\t$, then
\[
\lb{k0} {|f(k)|\/|2k|}\ge {6-\sqrt e\/8}>{1\/2}.
\]
Moreover,  the number of zeros $\cN(\r)$ of $f$ (counted with
multiplicity) in disk $\dD_\r(it)$ with the center $it=i2\t$ and the
radius $\r\ge \sqrt 2  \t $ satisfies
\[
\lb{zf} \cN(\r) \le 1+{4\/\log 2}\rt({\g \r\/\pi}+{\t\/\r}\rt),
\]
In particular, the number of zeros $\cN_+$ of $f$ (counted with
multiplicity) in $\ol \C_+$  satisfies
\[
\lb{qw2} \cN_+ \le C_1+C_2 \g \t,
\]
where  the constants $C_1=1+{8\/\sqrt{17}\log 2}+\ve_1$ and
$C_2={2\sqrt{17}\/\pi\log 2}+\ve_2, $ for any small $\ve_1,
\ve_2>0$.

\end{lemma}

\no {\bf Proof.} i) We have $\vs(k)\le {1\/2}$ and $\p={f\/2ik}$
satisfies $|\p(k)-(1-{\t_o\/2ik})|\le {\t \vs\/2|k|} e^{\vs}$ and
then $|\p(it)-(1+{\t_o\/2t})|\le {1\/8} \sqrt e$ for any $k\in
\ol\C_+$, which yields $|\p(k)|\ge 1- {|\t_o|\/4\t}-{1\/8} \sqrt e
\ge {6-\sqrt e\/8}>{1\/2}$.

ii)  Recall the Jensen formula (see p. 2 in \cite{Koo88}) for an
entire function $F$ any $r>0$:
\[
\lb{qw3} \log |F(0)|+\int _0^r{\cN_s (F)\/s}ds= {1\/ 2\pi }\int
_0^{2\pi}\log |F(re^{i\f})|d\f,
\]
where $ \cN_s (F)$ is the number of zeros of  $F$ in the disk
$\dD_s(0)$. We take the function $F(z)=f(it+z)$ and the disk
$\dD_r(it), t=2\t$ with the radius $r\ge \sqrt 8\t$. Thus \er{qw3}
implies
\[
\begin{aligned}
\lb{qw5} \log |f(it)|+\int _0^r{\cN (s)\/ s}ds={1\/ 2\pi }\int
_0^{2\pi}\log \big|f(k_\f)\big|d\f,\qq k_\f=it-ire^{i\f},
\end{aligned}
\]
where $\cN (s)=\cN_s (f(it+\cdot))$. We rewrite \er{qw5} in terms of
a function $\p(k):={f(k)\/2ik}$:
\[
\begin{aligned}
\lb{qw5x} \log |\p(it)|+\log |2t|+\int _0^r{\cN (s)\/ s}ds =S+ \log
{|2r|} , \qqq S={1\/ 2\pi }\int _0^{2\pi}\log \big|\p(k_\f)\big|d\f,
\end{aligned}
\]
where we have used the identity
$$
\int _\T\log {|it-re^{i\f}|}d\f=\int _\T\log
\rt(r\rt|{it\/r}-e^{i\f}\rt|\rt)d\f= 2\pi\log r,
$$
since $\int _\T\log |\a- e^{i\f}|d\f=0$ for any $|\a|<1$.
 From \er{dE} we obtain at $k_\f=it-ire^{i\f}$ and
 $a\in [{\pi\/4},{\pi\/2}]$ defined by $\cos a={t\/r} $:
\[
\lb{qw7}
\begin{aligned}
S={1\/ 2\pi }\int _0^{2\pi}\log \big|\p(k_\f)\big|d\f=S_0+S_1, \qqq
S_0={1\/ 2\pi }\int _{ -a}^{a}\log \big|\p(k_\f)\big|d\f,
\\
S_1={1\/ 2\pi }\int _a^{2\pi -a}\log \big|\p(k_\f)\big|d\f\le {1\/
2\pi }\int _a^{2\pi-a}\log \big(1+\vs_1 e^{\vs_1}\big)d\f={\pi-a
\/\pi}\log \big(1+\vs_1 e^{\vs_1}\big),
\end{aligned}
\]
where on the circle $\{z:|it-z|=r\}$ using \er{dE} we get for
$k_\f=it-ire^{i\f}$:
\[
\lb{qw10}
\begin{aligned}
& \min_{a\le\f\le \pi} |k_\f|=|t-re^{ia}|=r\sin a=\sqrt{r^2-t^2},
\\
& \vs_1:= \max_{a\le\f\le \pi}\vs(k_\f)={\t\/r\sin a}\le
{\t\/r}\sqrt2,
\end{aligned}
\]
and
\[
\lb{qw11}
\begin{aligned}
& \min_{0\le\f\le a} |k_\f|=(r-t)=r(1-\cos a),
\\
& \vs_o:= \max_{0\le\f\le a}\vs(k_\f)={\t\/r(1-\cos a)}\le
{\t\/r(1-{1\/\sqrt2})}.
\end{aligned}
\]
Consider the integral $S_{00}={2\g\/ \pi }\int _{0}^{a} k_- d\f$,
where we have $k_-=r\cos \f -t=r(\cos \f-\cos a)$ for $\f\in
(-a,a)$. Thus we obtain
\[
\lb{qw9}
\begin{aligned}
S_{00}= {2\g\/ \pi }\int _{0}^{a} r(\cos \f-\cos a) d\f={2 \g r\/
\pi } (\sin a-a\cos a)\le {2 \g r\/ \pi } ,
\end{aligned}
\]
and \er{dE}, \er{qw10} give
\[
\lb{qw8}
\begin{aligned}
S_0={1\/ 2\pi }\int _{ -a}^{a}\log \big|\p(k_\f)\big|d\f \le {1\/
\pi }\int _{0}^{a}\log \big(1+\vs_o e^{2\g k_-+\vs_o }\big)d\f
\\
\le S_{00}+ {1\/ \pi }\int _{0}^{a} \log \big(1+\vs_oe^{ \vs_o
}\big)\f=S_{00}+{2a \/\pi}\vs_o.
\end{aligned}
\]
Collecting \er{qw7}-\er{qw8} we obtain
\[
\lb{qw11z} S\le {2\/\pi}(\g  r+(\pi-a)\vs_1+a\vs_o).
\]
We  show \er{zf} and let $\r={r\/2}$. We have
$$
\int _0^r\cN (s){ds\/ s}\geq \cN (\r)\int _{\r}^r{ds\/ s}=\cN
(\r)\log 2.
$$
Summing \er{qw7}-\er{qw11z} we obtain
$$
-\log 2+\cN (\r)\log 2<{2\/\pi}(\g  r+(\pi-a)\vs_1+a\vs_o)\le
{2\/\pi}\g  r+{3\/2}\vs_1+{1\/2}\vs_o\le {2\/\pi}\g r+{8\t\/r},
$$
which yields \er{zf}. We  show \er{qw2}. We take any $r=\sqrt{17}
\t$ and $\r={r\/2}>{\sqrt{17}\/2}\t$. Then all zeros of $f$ in
$\ol\C_+$ belong to the disk $\ol\dD_{\r}(it)$ and  \er{zf} implies
$$
\cN_+\le C_1+\g \t C_2,\qq C_1=1+{8\/\sqrt{17}\log 2}+\ve_1, \qq
C_2={2\sqrt{17}\/\pi\log 2}+\ve_2,
$$
for any small $\ve_1, \ve_2>0$, which yields \er{qw2}.
 \BBox

 \

 \setlength{\unitlength}{1.0mm}
\begin{figure}[h]
\centering
\unitlength 1.0mm 
\begin{picture}(100,80)
\put(0,20){\vector(1,0){90.00}} \put(40,0){\vector(0,1){85.00}}
\put(31.0,82){$\Im k$} \put(83,16){$\Re k$}

\bezier{600}(40,72)(67,71)(68,44) \bezier{600}(40,72)(13,71)(12,44)
\bezier{600}(40,16)(67,17)(68,44) \bezier{600}(40,16)(13,17)(12,44)

\bezier{600}(40,76)(70.6,74.6)(72,44)
\bezier{600}(40,76)(9.4,74.6)(8,44)
\bezier{600}(40,12)(70.6,13.4)(72,44)
\bezier{600}(40,12)(9.4,13.4)(8,44)

\put(40,44){\vector(-1,0){28.00}}

\bezier{600}(40,31)(50,30)(51,20) \bezier{600}(40,31)(30,30)(29,20)

\put(40,44){\circle*{1}} \put(40,31){\circle*{1}}
\put(51,20){\circle*{1}} \put(29,20){\circle*{1}}
\put(62.5,20){\circle*{1}} \put(17.5,20){\circle*{1}}

\bezier{300}(40,44)(51.25,32)(62.5,20)
\bezier{300}(40,44)(28.75,32)(17.5,20)

\put(40.6,16.8){$0$} \put(62,16){$2Q$} \put(10,16){$-2Q$}
\put(40.5,32){$iQ$} \put(41.0,44){$t=2iQ$} \put(53.0,31.5){$r$}
\put(23.0,45){$r_1$} \put(2,3)

\linethickness{0.8pt}
\end{picture}
\caption{\footnotesize The case $2Q<r_1<2\sqrt{2}\,Q$} \lb{Fig1}
\end{figure}
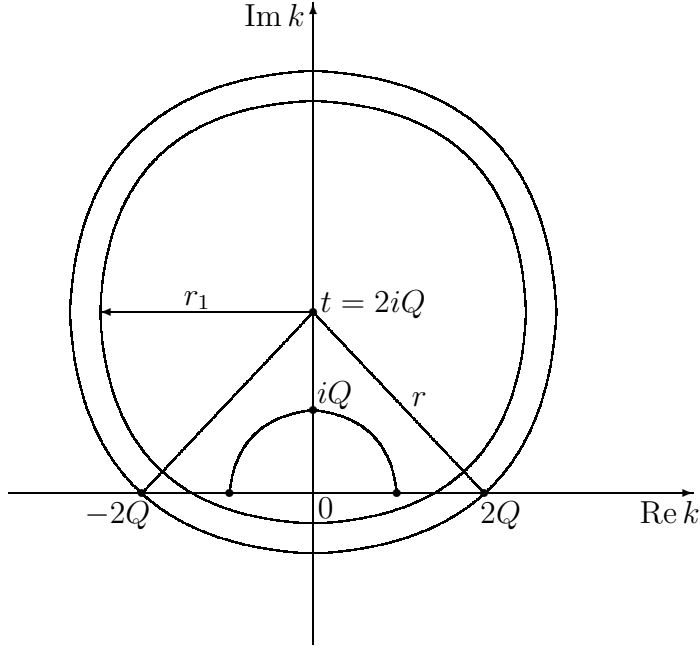

\subsection{Fundamental solutions.}
Consider  the Schr\"odinger operator $H$, when the potential $q\in
L^1(\R)$  is  complex and  $\supp q\ss [0,\g]$ for some $\g>0$. In
this case the equation \er{2.2} for $ y_+(x,k)=e^{-ikx}f_+(x,k) $
for all  $(x,k)\in [0,\g]\ts \C$ has the form
\[
\lb{y1} y_+(x,k)=1+\int _x^\g {\sin k
(t-x)\/k}e^{ik(t-x)}q(t)y_+(t,k)dt.
\]
 We recall well-known some properties of the functions $f_+, y_+$
(see e.g., \cite{K16}).

\begin{lemma}
\label{Tw1}
 Let $q\in L^1(\R)$ and  $\supp q\ss [0,\g ]$. Then the  function
 $f_+(0,k)$  is entire and satisfies
\[
\begin{aligned}
\lb{y}
|y_+(x,k)|\le e^{2(\g -x)k_-+\vs},\\
|y_+(x,k)-1|\le \vs e^{2(\g-x)k_-+\vs},
\end{aligned}
\]
for all $k\in \C$,  where $k_-={1\/2}(|\Im k| -\Im k)$ and $
\vs=\min \{\|q\|_1,{\|q\|\/|k|}\}$.

\end{lemma}

\begin{lemma}
\label{Tw2} Let $q\in L^1(\R)$ and  $\supp q\ss [0,\g ]$.  Then the
Wronskian $w$ is entire and satisfies
\[
\lb{w2z}
\begin{aligned}
& w(k)=2ik-q_0-\int_0^{\iy}q(t)(y_+(t,k)-1)dt,\\
\end{aligned}
\]
\[
\begin{aligned}
\lb{w3z}  |w(k)-(2ik-q_0)|\le \vs\|q\|e^{2\g  k_-+\vs(k)},
\end{aligned}
\]
where $q_0=\int_\R qdx$ and $\vs(k)=\min\{\|q\|_1,{\|q\|\/|k|}\}$.

\end{lemma}
{\bf Proof.} We show \er{w1}. From  Lemmas \ref{TL1}, \ref{TL2} we
have
\[
\lb{wx}
\begin{aligned}
w=\{f_-,f_+\}={f_+}'(0,k)+ikf_+(0,k)=2ik+\int_0^\g(i\sin kx-\cos
kx)q(x)f_+(x,k)dx,
\\
w-2ik=-\int_0^\g e^{-ikx}q(x)f_+(x,k)dx=-\int_0^\g q(x)y_+(x,k)dx.
\end{aligned}
\]
Thus using \er{wx} and Lemma \ref{Tw1} we obtain \er{w2z}, \er{w3z}.
\BBox

\subsection{Proof of Theorem \ref{Tcp}}
Consider  the Schr\"odinger operator $H$, when the potential $q\in
L^1(\R)$  is  complex and
 $\supp q\ss [0,\g]$ for some $\g>0$. Due to Lemma \ref{Tw2} the Wronskian
 $w$  is entire and satisfies
$
 |w(k)-(2ik-q_0)|\le \vs(k) e^{2\g k_- +\vs(k)}
$ for all $k\in \C$, where $\vs(k)=\min \{\|q\|_1, {\|q\|\/|k|_1}\}$
  Thus due to this fact and \er{1.2} we obtain that the function
$w\in \cE_\g$ with $\t=\|q\|, \t_o=q_0$ and $\s=\|q\|_1$. Then Lemma
\ref{Tef} gives the estimate \er{com1}, \er{com2}.
 \BBox


\section {Operators on a half-line with the Neumann boundary condition}
\setcounter{equation}{0}

\subsection{Definitions}

In this section we consider Schr\"odinger operators $H_{\bf {_n}}$
on $L^2(\R_+)$ and $H_{\bf {_d}}$  given by
$$
\begin{aligned}
 & Case \ 1: \qqq H_{\bf
{_n}}y=-y''+q_{{}_+}y,\qq
{\rm with \ Neumann \ b.c.} \ y'(0)=0,\\
& Case \ 2: \qqq  H_{\bf {_d}}y=-y''+q_{{}_+}y,\qq  {\rm with \
Dirichlet \ b.c.} \ y(0)=0,
\end{aligned}
$$
where the potential $q_{{}_+}$ is complex and satisfies:
\[
\lb{NdV} \int_0^\iy (1+x)|q_{{}_+}(x)|dx<\iy.
\]
We  define functions
 $w_{\bf {_n}}, \p_{\bf {_n}}, \P_{\bf {_n}}$ and
$w_{\bf {_d}}=\p_{\bf {_d}}$ in $\C_+$ by
$$
 w_{\bf {_n}}
=f_+'(0,k),\qq \p_{\bf {_n}}(k)={w_{\bf {_n}}\/ik},\qq \P_{\bf
{_n}}={w_{\bf {_n}}\/i(k+i)},\qqq \p_{\bf {_d}}=w_{\bf {_d}}
=f_+(0,k)
$$
 It is known that the operator $H_a, a={{ \bf
n, d}}$ has two components: the essential spectrum which covers the
half-line $[0,\iy)$
 plus $N_a \le \iy$ eigenvalues (counted with multiplicity) in the cut
 spectral domain $\C\sm [0,\iy)$. We denote them by
 $
 E_{a,j}\in \C\sm [0,\iy), j=1,...,N_a,
 $
according to their multiplicity. Note, that the multiplicity of each
eigenvalue equals 1, but we call the multiplicity of the eigenvalue
its algebraic multiplicity.  We call $k_{a,j}=\sqrt E_{a,j}\in \C_+$
also the eigenvalues of the operator $H_a$ and label them by $ \Im
k_{a,1}\ge \Im k_{a,2}\ge \Im k_{a,3}\ge... $. We define the set
$\dk_a=\{k_{a,1},k_{a,2},...\in \C_+\}$ and the Blaschke product
$$
B_a=\prod_{j\ge 1}\frac{k-k_{a,j}}{ k-\ol k_{a,j}},\qqq k\in \C_+.
$$
 This product converges absolutely for each $k\in
\C_+$, since all zeros of $w_a$ are uniformly bounded, see \er{df1}.
 Moreover, it has an analytic continuation from
$\C_+$ into the domain $\{|k|>r_c\}$, where $r_c={\|q\|\/2}$ and has
the following Taylor series
\[
\begin{aligned}
\lb{NB6} & \log
B_a(k)=-i{B_{a,0}\/k}-i{B_{a,1}\/2k^2}-i{B_{a,2}\/3k^3}-..., \qqq as
\qqq |k|>r_c,
\end{aligned}
\]
where $B_{a,0}=2\sum_{j=1}^N\Im k_{a,j},....$

We show the relations between the operators $H$ and $H_a$. Here we
use a standard trick. For $q_{_+}\in L^1(\R_+)$
 we define the Schr\"odinger operator $\wt H$ on $L^2(\R )$
with an even potential $\wt q$ by
\[
\lb{Nwtq} \wt Hy=-y''+\wt qy,\qqq \qqq \wt q(x)=q_{_+}(|x|), x\in
\R.
\]
 Recall that the Schr\"odinger equation $
  -f''+\wt qf=k^2f,
$ has unique Jost solutions $ \wt f_{\pm}(x,k)$. For the case of the
even potential $\wt q$ they satisfy
$$
\wt f_+(x,k)=\wt f_-(-x,k) \qqq \forall x\ge 0.
$$
This implies that the Wronskian $\wt w(k)$ for the operator  $\wt H$
has the specific form:
\[
\lb{Ntw} \wt w(k)=\{\wt f_+(x,k), \wt f_-(x,k)\}|_{x=0}=2 f_+(0,k)
f_+'(0,k)=2 w_{\bf {_d}}w_{\bf {_n}}.
\]

 This identity is
very useful. For example, if we have estimate \er{com1}, \er{com2}
for the operator $H$, then \er{Ntw} gives the estimate \er{com1x},
\er{com2x} for the operators $H_{\bf {_n}}$ and $H_{\bf {_d}}$.
Consider  estimates for  complex compactly supported potentials
$q_{_+}$. In this case the Jost function $\p(k)=f_+(0,k )$ is entire
and due to \er{1.2} it has a finite number of zeros in $\C_+$.

  \begin{corollary} \lb{TNcp}
Let $q_{{}_+}\in L^1(\R_+)$ and let $\supp q_{{}_+}\ss [0,\g]$ for
some $\g>0$.    Then the number of zeros $\cN_\r(w_a)$ of $w_a,
a={\bf n, d}$ (counted with multiplicity) in disk $\dD_\r(it)$ with
the center $it=i4\|q_{_+}\|$ and the radius $\r\ge \sqrt 8
\|q_{_+}\| $ satisfies
\[
\lb{com1x} \cN_\r(w_{\bf {_n}})+\cN_\r(w_{\bf {_d}}) \le 1+{4\/\log
2}\rt({\g \r\/\pi}+{2\|q_{_+}\|\/\r}\rt).
\]
In particular, the number of zeros $\cN_+(w_a)$ of $w_a$ (counted
with multiplicity) in $\ol\C_+$ satisfies
\[
\lb{com2x} \cN_+(w_{\bf {_n}})+\cN_+(w_{\bf {_d}})     \le C_1
+2C_2\g \|q_{_+}\|,
\]
where the constants $C_1, C_2\le 5$ (see more about $C_1, C_2$ in
Lemma \ref{Tef}).

\end{corollary}

 {\bf Proof.} Note that we have
$ \supp \wt q\ss [-\g,\g]$ and the identity $\|\wt q\|=2\|q_+\|. $
 Then using  these relations and \er{Ntw} and \er{com1}, \er{com2}
for the operator $H$, we obtain \er{com1x}, \er{com2x} for the
operators $H_{\bf {_n}}$ and $H_{\bf {_d}}$. \BBox


 We describe the basic properties of $H_{\bf {_d}}, H_{\bf {_n}}$
and the Blaschke product $B_{\bf {_n}}$.


\begin{proposition}\lb{TN1}
Let a  potential $q_{_+}$ be complex and let $c_0:=\int_{\R_+}
|q_{_+}(x)|dx<\iy$. Then

\no i) Each eigenvalue $E$ of $H_{\bf {_n}}$ satisfies
$|E|^{1\/2}\le c_0$.

\no ii) Let $q_{{{}_+}, 0}=\int_0^\iy q_{_+}(x)dx$ and $
c_1:=\int_{\R_+} |xq_{_+}(x)|dx<\iy$,  and let $A=c_0c_1e^{c_1}$.
Then

$\bu$ If $A<\Re q_{{{}_+}, 0}$, then the operator $H_a$ does  not
have eigenvalues.

$\bu$  If $A<\Re (-q_{{{}_+}, 0})$, then   $H_{\bf {_d}}$ does not
have eigenvalues and the operator $H_{\bf {_n}}$ has one simple
eigenvalue.

\no iii) The Blaschke product $B_{\bf {_n}}(k), k\in \C_+$ given by
\er{Bk} belongs to $\mH_\iy$ with $\|B_{\bf {_n}}\|_{\mH_\iy}\le 1$.

\end{proposition}

 {\bf Proof.} The proof is similar to the case of Proposition
\ref{T1}. \BBox

\subsection{Trace formulas and estimates}
We describe the Jost function $\p_{\bf {_n}}$ in terms of a
canonical factorization, which in  general $\p_{\bf {_n}}\notin
\mH_\iy$.

\begin{theorem}\label{TN2}
Let $q_{_+}$ satisfy \er{NdV}. Then $\p_{\bf {_n}}$ has a standard
canonical factorization:
$$
\p_{\bf {_n}} =\p_{{\bf {_n}},in} \p_{{\bf {_n}},out},
$$
where the inner factor $\p_{{\bf {_n}},in}$ and the outer factor
$\p_{{\bf {_n}},out}$ are given by
\[
\lb{NDi1}
\begin{aligned}
 \p_{{\bf {_n}},in}(k)={B_{\bf {_n}}(k)\/ik}e^{-iK_{\bf {_n}}(k)},\qqq  K_{\bf {_n}}(k)={1\/\pi}\int_\R
{d\n_{\bf {_n}}(t)\/k-t}, \\
\p_{{\bf {_n}},out}(k)=e^{iM_a(k)},\qqq  M_{\bf {_n}}(k)=
{1\/\pi}\int_\R {\log |\p_{\bf {_n}}(t+i0)|\/k-t} dt,
\end{aligned}
\]
for all $k\in \C_+$ and the function $\log |\p_{\bf {_n}}(t+i0)|$
belongs to $L_{loc}^1(\R)$.

\no $\bu$  $d\n_{\bf {_n}}(t)\ge 0$ is some singular compactly
supported measure on $\R$, which satisfies
\[
\lb{Nsmsx}
\begin{aligned}
\n_{\bf {_n}}(\R)=\int_\R d\n_{\bf {_n}}(t)<\iy,\qqq
 \supp \n_{\bf {_n}} \ss \{z\in \R: w_{\bf {_n}}(z)=0\}\ss [-r_c, r_c],
 \end{aligned}
\]
$\bu$  The function $K_{\bf {_n}}(\cdot)$ has an analytic
continuation from $\C_+$ into the cut domain $\C\sm [-r_c, r_c]$
and has the following Taylor series in the domain $\{|k|>r_c\}$:
\[
\lb{NKn}
\begin{aligned}
 K_{\bf {_n}}(k)=\sum_{j=0}^\iy {K_{{\bf {_n}},j}\/k^{j+1}},\qqq  \qq K_{{\bf {_n}},j}={1\/\pi}\int_\R
t^jd\n_{\bf {_n}}(t).
 \end{aligned}
\]
\end{theorem}

\no {\bf Proof.} The proof is similar to the case of Theorem
\ref{T2}.
 \BBox

The Jost function $\p_{\bf {_d}}=f_+(0,k)$ for the operator $H_{\bf
{_d}}$ on $\R_+$ with the Dirichlet boundary condition at $x=0$ also
the has a standard canonical factorization similar to $\p_{\bf
{_n}}$ in Theorem \ref{TN2}, see \cite{K18}. Moreover, the
corresponding inner factor $\p_{{\bf {_d}},in}$ is expressed in
terms of some singular compactly supported measure $d\n_{\bf
{_d}}(t)\ge 0$ on $\R$, which satisfies
\[
\lb{Nsd}
\begin{aligned}
\n_{\bf {_d}}(\R)=\int_\R d\n_{\bf {_d}}(t)<\iy,\qqq
 \supp \n_{\bf {_d}} \ss \{z\in \R: \p_{\bf {_d}}(z)=0\}\ss [-r_c, r_c],
 \end{aligned}
\]
Then from \er{Nsmsx}  and \er{Nsd} we obtain the simple fact:
\[
\lb{Nss}
\begin{aligned}
\supp \n_{\bf {_n}}\cap \supp \n_{\bf {_d}}=\es.
 \end{aligned}
\]

\medskip

  \begin{theorem} \lb{TN3}
Let a potential $q_{_+}$ satisfy \er{NdV}. Then
\[
\begin{aligned}
\lb{Ntr1} B_{{\bf {_n}},0}+{\n_{\bf {_n}}(\R)\/\pi}+{1\/2}\int_0^\iy
\Re q_{{}_+}(x)dx={\rm v.p.}{1\/\pi}\int_\R \log |\p_{\bf {_n}}(t)|
dt,
\end{aligned}
\]
where the integral in the r.h.s converges. Moreover, the following
hold true
\[
\begin{aligned}
\lb{NeBs} B_{{\bf {_n}},0}+{\n_{\bf {_n}}(\R)\/\pi}+{1\/2}\int_0^\iy
\Re q_{{}_+}(x)dx\le \ {2\/\pi}\big(1+\|q_+\|_1\big)
+r_+\big(C_o+\log r_+\big),
\end{aligned}
\]
where $\|q_+\|_1=\int_0^\iy |xq_+(x)|dx$ and $ r_+={\|q_+\|}$ and
$C_o$ is some absolute constant.

  \end{theorem}

  \no {\bf Proof.} The proof is similar to the case of Theorem
\ref{T3} and \ref{T4}.  \BBox

\medskip

\medskip

\footnotesize\footnotesize

\no {\bf Acknowledgments.} \footnotesize  EK is grateful to Ari
Laptev (London) for discussions about the Schr\"odinger operators
with complex potentials and to Alexei Alexandrov (St.Petersburg) for
discussions and useful comments about Hardy spaces. Our study was
supported by the RSF grant No 18-11-00032.

\end{document}